\documentclass[a4paper,leqno,11pt,final]{article}
\usepackage{latexsym}
\usepackage{amsthm,amsmath,amssymb}

\usepackage[dvips]{graphicx, color}

\newtheorem{thm}{Theorem}

\newtheorem{prop}[thm]{Proposition}

\newtheorem{conj}[thm]{Conjecture}
\newtheorem{lem}[thm]{Lemma} 

\newcommand{\Aut}{\mathrm{Aut}}

\newcommand{\het}{\mathrm{ht}}
\newcommand{\id}{\mathrm{id}}

\newcommand{\Res}{\mathrm{Res}}

\newcommand{\Vol}{\mathrm{Vol}}
\newcommand{\W}{\frak{W}}

\begin{document}
\title{\bf {\Large{Zeta Functions  for Function Fields}}}
\author{\bf Lin WENG}
\date{\bf }
\maketitle
\noindent
{\bf Abstract.} {\footnotesize To (weighted) count semi-stable bundles
on curves defined over finite fields, we introduce new genuine zeta functions. 
There are two types, i.e., the pure non-
abelian zetas defined using  semi-stable bundles, and the group zetas 
defined for pairs consisting of (reductive group, maximal parabolic subgroup). 
Basic properties such as rationality and functional equation
are obtained. Moreover,  conjectures on their zeros and uniformity
are given.
\vskip 0.05cm
The constructions and results  were announced in our paper on \lq Counting Bundles'.}

\section{Pure Non-Abelian Zeta Functions}
Non-Abelian zeta functions for function fields  were introduced in [W1] about 10 years ago. 
However, due to the lack of the Riemann Hypothesis, we have faced some essential difficulties. 
Recently, with an old paper of Drinfeld ([D]) on counting rank two cuspidal  $\mathbb Q_l$-representations 
for function fields, we realize that our old 
definition of zeta should be altered: instead of counting rational semi-stable bundles of all degrees as done in [W], 
only these of degree zero and more generally degrees of multiples of the rank should be counted. 
This then leads to the definition of {\it pure} high rank zeta functions of this paper. This purity proves to be very essential: We expect that the Riemann Hypothesis holds for pure zetas. Indeed, in this direction, we now have the work of Yoshida [Y] for the RH in rank two and of mine [W4] for elliptic curves.

\subsection{Counting Semi-Stable Bundles}
Let $X$ be an irreducible, reduced, regular projective curve of genus $g$ defined over $\mathbb F_q$. Denote by $\mathcal M_{X,r}(d)$ the moduli space of rank $r$ semi-stable bundles of degree $d$ consisting of the Seshadri Jordan-H\"older equivalence classes of $\mathbb F_q$-rational semi-stable bundles on $X$. For our own purpose, we consider $\mathcal M_{X,r}(d)$ in the sense of the {\it fat} moduli, meaning that ordinary moduli spaces equipped with an additional structure at Seshadri class $[\mathcal E]$
defined by the collection of semi-stable bundles in $[\mathcal E]$, namely, the set $\big\{\mathcal E:\mathcal E\in [\mathcal E]\big\}$, is added at the point $[\mathcal E]$. For our own convenience,
$\mathcal M_{X,r}(d)$ eqiupped with such a structure is called a {\it fat moduli space} and denoted as ${\bf M}_{X,r}(d)$. 

A natural question is to count these 
$\mathbb F_q$-rational semi-stable bundles $\mathcal E$ on $X$. For this purpose, two naive invariants, namely, the automorphism group $\Aut(\mathcal E)$ and its global sections $h^0(X,\mathcal E)$, will be used.
This then leads to the refined Brill-Noether loci $$W_{X,r}^i(d):=\Big\{[\mathcal E]\in {\bf M}_{X,r}(d):
\min_{\mathcal E\in[\mathcal E]}\{h^0(X,\mathcal E)\}\geq i\Big\}$$ and $$[\mathcal E]^j:=\{\mathcal E\in[\mathcal E]: \dim_{\mathbb F_q}\Aut\,{\mathcal E}\geq j\}.$$ Recall that there exist natural
isomorphisms $${\bf M}_{X,r}(d)\to {\bf M}_{X,r}(d+rm),\qquad
\mathcal E\mapsto A^m\otimes\mathcal E$$ and $${\bf M}_{X,r}(d)\to {\bf M}_{X,r}(-d+r(2g-2)),\qquad
\mathcal E\mapsto K_X\otimes\mathcal E^\vee,$$ 
where $A$ is an Artin line bundle of degree one on $X/\mathbb F_q$ and $K_X$ denotes the dualizing bundle of $X/\mathbb F_q$. So we only need to count ${\bf M}_{X,r}(d_0)$ for $d_0= 0, 1,\dots, r(g-1).$ Accordingly, we introduce
$$\alpha_{X,r}(d):=\sum_{\mathcal E\in {\bf M}_{X,r}(d)}\frac{q^{h^0(X,\mathcal E)}-1}{\Aut(\mathcal E)},\qquad \beta_{X,r}(d):=\sum_{\mathcal E\in {\bf M}_{X,r}(d)}\frac{1}{\Aut(\mathcal E)}$$ with $\beta$ a classical invariant ([HN]).

So to count bundles, the problem becomes how to control $\alpha_{X,r}(d_0)$'s with $d_0$ ranging as above, and 
$\beta_{X,r}(d)$ with $d=0,1,\dots r-1$. For $\alpha$, two general principles can be used for counting semi-stable bundles, namely, 

\noindent
(i) the vanishing theorem claiming that, for semi-stable $\mathcal E$,
$$h^1(X,\mathcal E)=0\qquad\mathrm{if}\qquad d(\mathrm E)\geq r(2g-2)+1;$$

\noindent
(ii) the Clifford lemma claiming that, for semi-stable $\mathcal E$,
$$h^0(X,\mathcal E)\leq r+\frac{d}{2}\qquad\mathrm{if}\qquad 0\leq \mu(\mathcal E)\leq 2g-2.$$

By contrasting, the invariant $\beta$ has already been understood, thanks to the high profile works of Harder-Narasimhan ([HN]),
Desale-Ramanan ([DR]), Atiyah-Bott ([AB]), Witten ([Wi]) and Zagier ([Z]). 
To state it, let $$\zeta_X(s):=\frac{\prod_{i=1}^{2g}(1-\omega_iq^{-s})}{(1-q^{-s})(1-qq^{-s})}$$ be the Artin zeta function of $X/\mathbb F_q$,
 $$v_n(q):=\frac{\prod_{i=1}^{2g}(1-\omega_i)}{q-1}q^{(r^2-1)(g-1)}\zeta_X(2)\cdots\zeta_X(r)$$ and for a partition $r=n_1+\cdots+n_k$ of $r$, set
 $$c_{r,d}(t):=\prod_{i=1}^{s-1}\frac{t^{(n_i+n_{i+1})\{n_1+\cdots+n_i)d/n\}}}{1-t^{n_i+n_{i+1}}}.$$
  
 \begin{thm} ([HN], [DS], in particular, [Z, Thm 2]) For any pair $(r,d)$, we have
 $$\beta_{X,r}(d)=\sum_{n_1,\dots,n_s>0,\ \sum n_i=r}q^{(g-1)\sum_{i<j}n_in_j}c_{r,d}(q)\cdot \prod_{i=1}^s v_n(q).$$
 \end{thm}
\subsection{Pure Non-Abelian Zeta Functions}
Practically, the difficulty of counting semi-stable bundles comes form the fact that direct summands of the associated 
Jordan-H\"older graded bundle, or equivalently, the Jordan-H\"older filtrations, of an 
$\mathbb F_q$-rational semi-stable bundle in general would not be defined over $X/\mathbb F_q$, but
rather its scalar extension $X_n/\mathbb F_{q^n}$. Theoretically, this is the junction point where the abelian and non-abelian ingredients of curves interact. For examples, torsions of Jacobians, Weierastrass points and stable but not absolutely stable bundles are closely related and hence get into the picture naturally. Good examples may be found in [W1]. 

To uniformly study $\alpha$ and $\beta$'s, we, in [W], introduce the non-abelian zeta functions with the hope that
the Riemann Hypothesis would hold for them. Unfortunately, examples shows that there are zeros off the central line 
for these old zetas. (For details, see the examples below.) 
This, in practice, has prevented any further studies for such zetas. However, during my visit to IHES in September 2011,  
we got to know the work of Drinfeld ([D]). Learning from it, we now know where the problem lies for old zetas:
We should {\it  count only the pure part, instead of counting all}.
\vskip 0.30cm
\noindent
{\bf Main Definition 1.} {\it For an irreducible, reduced, regular projective curve $X$ of genus $g$ defined over finite field $\mathbb F_q$,  define its rank r pure non-abelian zeta function  by
$$\begin{aligned}\zeta_{X,r}(s):=&\sum_{m=0}^\infty
\sum_{V\in {\bf  M}_{X,r}(d), d=rm}\frac{q^{h^0(C,V)}-1}{\#\mathrm{Aut}(V)}\cdot (q^{-s})^{d(V)},\\\widehat\zeta_{X,r}(s):=&\sum_{m=0}^\infty
\sum_{V\in {\bf  M}_{X,r}(d), d=rm}\frac{q^{h^0(C,V)}-1}{\#\mathrm{Aut}(V)}\cdot (q^{-s})^{\chi(X,V)}.\end{aligned}$$ As usual, set
$$\qquad Z_{X,r}(t):=\zeta_{X,r}(s)\qquad\mathrm{and}\qquad \widehat Z_{X,r}(t):=\widehat\zeta_{X,r}(s)\qquad\mathrm{with}\qquad t:=q^{-s}.$$}

By [W1, Prop 1.2.1], when the rank is one,
$$\zeta_{X,1}(s)=\zeta_{X}(s)$$ is the Artin zeta function. 

Moreover
 $$\begin{aligned}\widehat Z_{X,r}(t)
=&\sum_{m=0}^{2g-2}
\sum_{V\in {\bf  M}_{X,r}(d), d=rm}\frac{q^{h^0(C,V)}-1}{\#\mathrm{Aut}(V)}\cdot t^{\chi(C,V)}\\
&+\sum_{m\,>\, 2g-2}
\sum_{V\in {\bf  M}_{X,r}(d), d=rm}\frac{q^{h^0(C,V)}-1}{\#\mathrm{Aut}(V)}\cdot t^{\chi(C,V)}\\
=&\Big(\sum_{m=0}^{(g-1)-1}
\sum_{V\in {\bf  M}_{X,r}(d), d=rm, r[(2g-2)-m]}\frac{q^{h^0(C,V)}-1}{\#\mathrm{Aut}(V)}\cdot t^{\chi(C,V)}\\
&\qquad+\sum_{V\in {\bf  M}_{X,r}(d), d=r(g-1)}\frac{q^{h^0(C,V)}-1}{\#\mathrm{Aut}(V)}\cdot t^{\chi(C,V)}\Big)\\
&+\sum_{m\,>\, 2g-2}
\sum_{V\in {\bf  M}_{X,r}(0)}\frac{q^{r[m-(g-1)]}-1}{\#\mathrm{Aut}(V)}\cdot t^{r[m-(g-1)]}\\
&\hskip 1.0cm\mathrm{(by\ the\ Vanishing\ Thm\ and\ the\ Riemann-Roch\ Thm)}\\
=&\Big[\sum_{m=0}^{(g-1)-1}
\Big(\sum_{V\in {\bf  M}_{X,r}(rm)}\frac{q^{h^0(C,V)}-1}{\#\mathrm{Aut}(V)}\cdot t^{r[m-(g-1)]}\\
&\qquad+\sum_{W\in {\bf  M}_{X,r}(rm)}\frac{q^{h^0(C,W)-r[m-(g-1)]}-1}{\#\mathrm{Aut}(W)}\cdot t^{r[(g-1)-m]}\Big)\\
&\qquad+\sum_{V\in {\bf  M}_{X,r}(r(g-1))}\frac{q^{h^0(C,V)}-1}{\#\mathrm{Aut}(V)}\cdot t^{0}\Big]\\&+\sum_{V\in {\bf  M}_{X,r}(0)}\frac{1}{\#\mathrm{Aut}(V)}\cdot\Big(\frac{(qt)^{rg}}{1-(qt)^r}-\frac{t^{rg}}{1-t^r}\Big)\\
&\hskip 7.5cm\mathrm{(by\ the\ duality)}\\
=&\Big[\sum_{m=0}^{(g-1)-1}
\Big(\sum_{V\in {\bf  M}_{X,r}(rm)}\frac{q^{h^0(C,V)}-1}{\#\mathrm{Aut}(V)}\cdot t^{r[m-(g-1)]}\\
&+\sum_{W\in {\bf  M}_{X,r}(rm)}\frac{q^{h^0(C,W)}-1}{\#\mathrm{Aut}(W)}\cdot ({qt})^{r[(g-1)-m]}
\\
&+\sum_{W\in {\bf  M}_{X,r}(rm)}\frac{q^{-r[m-(g-1)]}-1}{\#\mathrm{Aut}(W)}\cdot t^{r[(g-1)-m]}\Big)+\sum_{V\in {\bf  M}_{X,r}(r(g-1))}\frac{q^{h^0(C,V)}-1}{\#\mathrm{Aut}(V)}\cdot t^{0}\Big]\\
&+\sum_{V\in {\bf  M}_{X,r}(0)}\frac{1}{\#\mathrm{Aut}(V)}\cdot\Big(\frac{(qt)^{rg}}{1-(qt)^r}-\frac{t^{rg}}{1-t^r}\Big)\\
=&\Big[\sum_{m=0}^{(g-1)-1}
\alpha_{X,r}(rm)\cdot \Big(t^{r[m-(g-1)]}+ \Big(\frac{1}{qt}\Big)^{r[m-(g-1)]}\Big)+\alpha_{X,r}(r(g-1))\Big]\\
&\qquad+\beta_{X,r}(0)\cdot \Big(\frac{(qt)^r}{1-(qt)^r}-\frac{t^r}{1-t^r}\Big)
\end{aligned}$$
Consequently,
$$\widehat Z_{X,r}(\frac{1}{qt})=\widehat Z_{X,r}(t),$$
and, if we introduce $T:=t^r$ and $Q:=q^r$,
$$\begin{aligned}Z_{X,r}&(t)
=\sum_{m=0}^{(g-1)-1}\alpha_{X,r}(mr)
\cdot \Big(T^{m}+ Q^{(g-1)-m}\cdot T^{2(g-1)-m}\Big)\\
&+\alpha_{X,r}\big(r(g-1)\big)\cdot T^{g-1}+(Q-1)\beta_{X,r}(0)\cdot \frac{T^{g}}{(1-T)(1-QT)}.
\end{aligned}\eqno(1)$$
This then completes the proof of the following

\begin{thm} ({\bf Zeta Facts}) (i) $\zeta_{X,1}(s)=\zeta_X(s)$, the Artin zeta function for $X/\mathbb F_q$;

\noindent
(ii) ({\bf Rationality}) There exists a degree $2g$ polynomial $P_{X,r}(T)\in\mathbb Q[T]$ of $T$
such that
$$Z_{X,r}(t)=\frac{P_{X,r}(T)}{(1-T)(1-QT)}\quad\mathrm{with}\quad T=t^r,\ Q=q^r;$$
\noindent
(iii) ({\bf Functional equation}) $$\widehat Z_{X,r}(\frac{1}{qt})=\widehat Z_{X,r}(t),$$

\end{thm}
These non-abelian zetas give a systematical treatment of invariants $\alpha$'s and $\beta$'s in counting semi-stable bundles.
With $\beta_{X,r}(0)$ known, we expect the uniform control of $$\alpha_{X,r}(0), \alpha_{X,r}(0),\dots, \alpha_{X,r}\big(r(g-1)\big)$$ through the following

\noindent
{\bf Riemann Hypothesis.} {\it Let $P_{X,r}(T)=P_{X,r}(0)\cdot\prod_{i=1}^{2g}(1-\omega_{X,r}(i)T)$, then $$|\omega_{X,r}(i)|=Q^{\frac{1}{2}},\qquad\forall i,\ 1\leq i\leq 2g.$$}

\noindent
{\bf Examples.} (i) ([W4]) {\it Rank 2 Zeta for Elliptic Curves} Let $E$ be an elliptic curve defined over $\mathbb F_q$ with $N$ the number of $\mathbb F_q$-rational points. For rank two pure zeta, it suffices to calculate 
$\alpha_{E,2}(0)$ and $\beta_{E,2}(0)$.  By Thm 1, $$\beta_{E,2}(0)=\frac{N}{q-1}\Big(1+\frac{N}{q^2-1}\Big).$$
On the other hand, by the classification of Atiyah ([A]), over $\overline{\mathbb F_q}$, the graded bundle associated to a Jordan-H\"older filtration of a semi-stable bundle $V\otimes_{\mathbb F_q}\overline{\mathbb F_q}$ is of the form $\mathrm{Gr}(V\otimes_{\mathbb F_q}\overline{\mathbb F_q})=L_1\oplus L_2$ with $L_i$ degree zero line bundles, which may not be defined over $\mathbb F_q$. Consequently, for $\mathbb F_q$-rational semi-stable bundles $V$ of rank two, $h^0(E,V)\not=0$ if and only if $V=\mathcal O_E\oplus L$ or $V=I_2$ with $L$ a $\mathbb F_q$-rational line bundle of degree 0 and $I_2$ the only non-trivial extension of $\mathcal O_E$ by itself. Thus $\alpha_{E,2}(0)$ is given by $$\begin{aligned}&\Big(\frac{q^{h^0(E,\mathcal O_E\oplus \mathcal O_E)}-1}{\#\Aut(\mathcal O_E\oplus \mathcal O_E)}+\frac{q^{h^0(E,I_2)}-1}{\#\Aut (I_2)}\Big)+\sum_{L\in\mathrm{Pic}^0(E), L\not=\mathcal O_E}\frac{q^{h^0(E,\mathcal O_E\oplus L)}-1}{\#\Aut(\mathcal O_E\oplus L)}\\
=&\Big(\frac{q^2-1}{(q^2-1)(q^2-q)}+\frac{q-1}{(q-1)q}\Big)+(N-1)\frac{q-1}{(q-1)^2}=\frac{N}{q-1}.\end{aligned}$$ Thus, $$Z_{E,2}(t)=\alpha_{E,2}(0)\cdot\frac{1+(N-2)T+QT^2}{(1-T)(1-QT)},$$
the Riemann Hypothesis holds since $$\Delta=(N-2)^2-4Q=(N-2-2q)(N-2+2q)<0$$ using Hasse's theorem for the Riemann Hypothesis  of elliptic curves, namely $$N\leq 2\sqrt q.$$

\noindent
(ii) {\it Rank Two Bundles on Genus Two Curves} Let $X$ be a genus 2 curve. For rank two zeta, 
$$\begin{aligned}Z_{X,2}(t)
=&\alpha_{X,2}(0)
\cdot \Big(1+ Q\cdot T^{2}\Big)+\alpha_{X,r}\big(2\big)\cdot T\\
&+(Q-1)\beta_{X,2}(0)\cdot \frac{T^{2}}{(1-T)(1-QT)}.
\end{aligned}$$
Thus $$\begin{aligned}P_{X,2}(t)=&\alpha_{C,2}(0)\Big(1+Q^2T^4\Big)+\Big(\alpha_{X,2}(2)-
\alpha_{X,2}(0)(Q+1)\Big)\Big(T+QT^3\Big)\\
&+
\Big(2Q\,\alpha_{X,2}(0)-(Q+1)\alpha_{X,2}(2)+\beta_{X,2}(0)(Q-1)\Big)T^2.\end{aligned}$$
So we need to consider the spaces ${\bf M}_{X,2}(0)$, ${\bf M}_{X,2}(2)$ (and ${\bf M}_{X,2}(4)$). By the Clifford lemma,
$$h^0(X,V)=\begin{cases}0,1,2,&$if$\ \ V\in {\bf M}_{X,2}(0);\\
0,1,2,& $if$\ \ V\in {\bf M}_{X,2}(2);\\
2,3,4,&$if$\ \ V\in {\bf M}_{X,2}(4).\end{cases}$$ Consequently,
$$\begin{aligned}\alpha_{X,2}(0)=&\sum_{V\in W^1_{X,2}(0)}\frac{q-1}{\#\Aut(V)}+\frac{q^2-1}{(q^2-1)(q^2-q)},\\
\alpha_{X,2}(2)=&\sum_{V\in W^1_{X,2}(0)}\frac{q-1}{\#\Aut(V)}+\sum_{V\in W^2_{X,2}(0)}\frac{q^2-1}{\#\Aut(V)},\end{aligned}$$ note that $h^0(X,V)=2$ and $d(V)=0$ iff $V=\mathcal O_X\oplus\mathcal O_X$. Moreover, the Riemann Hypothesis now is equivalent to the  conditions that
$$A^2<4Q, \qquad B^2<4Q$$ with real constants $A,\ B$ defined by
 $$P_{X,2}(T)=\alpha_{X,2}(0)\cdot (1-AT+QT^2)(1-BT+QT^2).$$
That is to say, $$\begin{aligned}A+B=&(Q+1)-\alpha_{X,2}'(2),\\
AB=&(Q-1)\beta_{X,2}'(0)-(Q+1)\alpha_{X,2}'(2),\end{aligned}$$
where $$\alpha_{X,r}'(d):=\frac{\alpha_{X,r}(d)}{\alpha_{X,r}(0)},\qquad \beta_{X,r}'(0):=\frac{\beta_{X,r}(0)}{\alpha_{X,r}(0)}.$$

While the above does give a good control of $\alpha'$'s and $\beta'$, it looks a bit clumsy.
A much better way is to set $$Z_{X,r}(t)=\alpha_{X,r}(0)\cdot\exp\Big(\sum_{m=1}^\infty N_{X,r}(m)\frac{T^m}{m}\Big).$$
Then $$N_{X,r}(m)=1+Q^m-\sum_{i=1}^{2g}\omega_{X,r}(i)^m$$
and the Riemann Hypothesis gives a much elegant control of $N_{X,r}(m)$'s.
We expect that $N_{X,r}(m)$'s measure rank $r$ stable bundles over $X/\mathbb F_{q^m}$. This is certainly the case in rank one through Weil's counting zeta,
 and in rank two for  elliptic curves, as indicated in the above example.

\subsection{Why Purity} 

Next, we explain why purity is introduced for our study of zeta functions. Simply put,
this is due to the Riemann Hypothesis.

For this purpose, let $E$ be an irreducible, reduced regular elliptic curve defined over $\mathbb F_q$. We will concentrate on ranks two and three.  Due to isomorphisms $${\bf M}_{E,r}(d)\to {\bf M}_{E,r}(d+rm),\qquad
\mathcal E\mapsto A^m\otimes\mathcal E$$ and $${\bf M}_{E,r}(d)\to {\bf M}_{E,r}(-d),\qquad
\mathcal E\mapsto K_E\otimes\mathcal E^\vee,$$
among all invariants $\alpha$ and $\beta$'s, for rank two and three, 
it suffices to understand $\alpha_{E,r}(d)$ and $\beta_{E,r}(d)$ for $d=0,1$. 

\subsubsection{Rank Two}

From Ex(i) in \S1.2,
$$\beta_{E,2}(0)=\frac{N}{q-1}\Big(1+\frac{N}{q^2-1}\Big)\qquad\mathrm{and}\quad 
 \alpha_{E,2}(0)=\frac{N}{q-1}.$$
On the other hand, $\alpha_{E,2}(1)$ and $\beta_{E,2}(1)$ are easy to calculate.
Indeed, all semi-stable bundles of rank 2 and degree 1 are stable. Thus, by the classification of Atiyah ([A]),
${\bf M}_{E,2}(1)$ via the determinant line bundle map is isomorphic to $\mathrm{Pic}^1(E)$.
So $$\beta_{E,2}(1)=\frac{N}{q-1}.$$ Moreover, by the vanishing theorem, $h^0(E,\mathcal E)=1$.
Thus $$\alpha_{E,2}(1)=N.$$
As such, from [W1, \S1.2.2], we know that  
the original zeta function $\zeta$ of counting  all degree semi-stable bundles defined as
$$\sum_{V\in{\bf M}_{E,2}(d)}\frac{q^{h^0(E,V)}-1}{\#\Aut\,V}=:\zeta_{E,2}(s)+\zeta_{E,2}^1(s)$$
is given by
$$\begin{aligned}\Big(\alpha_{E,2}(0)+&\beta_{E,2}(0)\cdot\frac{(q^2-1)t^2}{(1-t^2)(1-q^2t^2)}\Big)+\beta_{E,2}(1)\Big(\frac{qt}{1-q^2t^2}-\frac{t}{1-t^2}\Big)\\
=&\frac{N}{q-1}\cdot\frac{1+(q-1)t+(N-1)t^2+(q-1)qt^3+q^2t^4}{(1-t^2)(1-q^2t^2)}.\end{aligned}$$
Here $$\zeta_{E,2}^1(s):=\sum_{V\in{\bf M}_{E,2}(2m+1),\ m\geq 0}\frac{q^{h^0(E,V)}-1}{\#\Aut\,V}.$$
Note that for the polynomial appeared in the numerator
$$P_2(t):=1+(q-1)t+(N-1)t^2+(q-1)qt^3+q^2t^4,$$ by the functional equation ([W1]), we have the 
factorization $$P_2(t)=(qt^2+A_+t+1)(qt^2+A_-t+1)$$ in $\mathbb R[t]$. Assume, as we may, that $|A_+|>|A_-|$. By Hasse's theorem for Artin zeta functions,
 the coefficients of $t^2$ in $P_2(t)$ is $N-1$, which is of the same order as 
$q-1$, the coefficient of $t$. Consequently,  $$A_+^2-4q>0,\qquad\mathrm{while}\qquad A_-^2-4q<0.$$
So there is no RH for the zeta defined by counting semi-stable bundles of all degrees.
\subsubsection{Rank Three}
We already saw that for pure rank two zeta functions, the RH holds. In fact, one can show that this patten persists ([W4]). This then leads to the problem of whether partial zeta functions defined by counting semi-stable bundles of other types of degrees satisfy the RH. Here, we use an example in rank 3 to indicate that another seemly natural choice does not work neither.

Introduce then the function
$$\zeta_{E,3}^{12}(s):=\sum_{\substack{V\in{\bf M}_{E,2}(d)\\
d\equiv 1,\,2\,(\mathrm{mod}\,3)}}\frac{q^{h^0(E,V)}-1}{\#\Aut\,V}.$$
(The reason for taking both 1 and 2, not just a single one, 1 or 2, in the congruence classes is that otherwise
the functional equation does not hold.)
By the vanishing theorem and the fact that all rank three semi-stable bundles of degree 1 or 2 are stable, 
one checks that
$$\zeta_{E,3}^{12}(s)=N\cdot \Big(\frac{qt+q^2t^2}{1-q^3t^3}-\frac{t+t^2}{1-t^3}\Big)=\frac{P_{E,3}(t)}{(1-t^3)(1-q^3t^3)}$$ where $$P_{E,3}^{12}(t)=(q-1)t\Big[q^2t^4+q(q-1)t^3+(q+1)t+1\Big].$$
The polynomial $$q^2t^4+q(q-1)t^3+(q+1)t+1$$ does not satisfy the Riemann Hypothesis.

\section{Group Zeta Functions}
 \subsection{Number Fields versus Function Fields}
For number fields, we have yet another type of zeta functions defined for pairs consisting of (reductive group, maximal parabolic subgroup)'s ([W2,3]). We will introduce such zeta functions for function fields next. 
For this purpose, we first examine analogue between function fields and number fields in our setting.
To be more precise, we will  analysis  
Zagier's formula for counting semi-stable bundles 
over curves on finite fields and our own volume formula for semi-stable lattices over number fields.

Set then $$\widehat\zeta_{F}(1):=\begin{cases}\Res_{s=1}\widehat\zeta_F(s), & F\ \mathrm{number\ field};\\
\Res_{s=1}\widehat\zeta_F(s)\cdot\log q, & F\ \mathrm{function\ field}.\end{cases}$$
 And denote by $\mathcal M_{\mathbb Q,r}[1]$ the moduli space of rank $r$
 semi-stable lattices of volume 1.
  
 \begin{thm} (i) (Reformulation of [Z, Thm 2]) For an irreducible, reduced, regular projective curve $X/\mathbb F_q$ of genus $g$,
 $$\frac{\beta_{X,r}(0)}{q^{-(g-1)\cdot\frac{r^2-r}{2}}}=\sum_{\substack{n_1,\dots,n_s>0,\\ n_1+\cdots+n_k=r}}
\frac{(-1)^{k-1}}{\prod_{j=1}^{k-1}(q^{n_j+n_{j+1}}-1)}\prod_{j=1}^k \prod_{i=1}^{n_j}\widehat\zeta_X(i);$$

\noindent
(ii) ([W2, \S4.8]) For a number field $F$, $$\frac{1}{r}\cdot \mathrm{Vol}\Big(\mathcal M_{\mathbb Q,r}[1]\Big)=\sum_{\substack{n_1,\dots,n_s>0,\\ n_1+\cdots+n_k=r}}
\frac{(-1)^{k-1}}{\prod_{j=1}^{k-1}(n_j+n_{j+1})}\prod_{j=1}^k \prod_{i=1}^{n_j}\widehat\zeta(i).$$
\end{thm}
Put Zagier's result in our form as above, the hidden parallel structures in these two worlds becomes crystal
 clear. That is to say, for the mass of moduli space of semi-stable objects, 
 when shift from number fields to function fields, the integers $n_j+n_{j+1}$ should be replaces by 
 $q^{(n_j+n_{j+1})}-1$. This then would suggest that, more generally, when defining group
 zeta functions associated to $(G,P)$ with $G$ reductive and  $P$ maximal parabolic for function fields, based on these for number fields investigated in [W2,3],
 we should replace the rational factor $\langle w\lambda-\rho,\alpha^\vee\rangle$ by
$q^{\langle w\lambda-\rho,\alpha^\vee\rangle}-1$.
In reality, even this is the direction we would go, this is not exactly the path we really pave. As a matter of fact, when shifting from number fields to function fields,  {\it the rational factor $\langle w\lambda-\rho,\alpha^\vee\rangle$ should be replaced by
$1-q^{-\langle w\lambda-\rho,\alpha^\vee\rangle}$, instead of $q^{\langle w\lambda-\rho,\alpha^\vee\rangle}-1$.}

\subsection{Definitions}
Let $X$ be an irreducible, reduced, regular projective curve of genus $g$ defined on $\mathbb F_q$. Denote by $F$ its function field. Let $G$ be a split connected reductive group with $B$ a fixed Borel over $F$.
Denote by $\Sigma(G):=\Sigma:=$
$$\Big(V,\langle\cdot,\cdot\rangle,\Delta=\{\alpha_1,\dots,\alpha_n\},\Lambda:=\{\lambda_1,\dots,\lambda_n\},\Phi=\Phi^+\cup\Phi^-,W\Big)$$ the associated root system with the Weyl vector  $\rho:=\frac{1}{2}\sum_{\alpha\in\Phi^+}\alpha$. For $w\in W$, set  
  $\Phi_w:=\Phi^+\,\cap\, w^{-1}\Phi^-$, and 
for $\alpha\in\Phi$, denote its coroot by $\alpha^\vee:=\frac{2}{\langle\alpha,\alpha\rangle}\cdot\alpha$.
 
From Lie theory, (see e.g., [H]), there is a well-known one-to-one correspondence between standard parabolic subgroups of $G$
and subsets of $\Delta$. Consequently, for a maximal standard parabolic subgroup $P$, there exists a unique $p=p(P)$ such that the  subset of $\Delta$ above for $P$ is given by $$\Delta_p=:\Delta\backslash\{\alpha_p\}:=\{\beta_{P,1},\dots,\beta_{P,n-1}\}.$$
For such $p$, let $\Phi_p$ be the corresponding root system. Then $\Phi_p$ is normal to the fundamental weight $\lambda_p$, since $\langle\alpha_i,\lambda_j\rangle=\delta_{ij}$. This similarly leads to positive $\Phi^+_p:=\Phi^+\,\cap\,\Phi_p$,  $\Phi_p^-,\,\rho_p$, and $W_p$.
Set $$c_p:=2\langle\lambda_p-\rho_p,\alpha_p^\vee\rangle.$$
Moreover, for $\lambda \in V$, introduce a specific coordinate system via $$\lambda=\sum_{j=1}^n(1+s_j)\lambda_j=\rho+\sum_{j=1}^ns_j\lambda_j.$$

\vskip 0.30cm
\noindent
{\bf Main Definition 2.}  {\it (i) The period of $G$ for $X$ is defined by 
$$\omega_X^G(\lambda):=\sum_{w\in W}\frac{1}{\prod_{\alpha\in\Delta}(1-q^{-\langle w\lambda-\rho,\alpha^\vee\rangle})}\prod_{\alpha\in\Phi_w}\frac{\widehat\zeta_X(\langle\lambda,\alpha^\vee\rangle)}
{\widehat\zeta_X(\langle\lambda,\alpha^\vee\rangle+1)}$$ where $\widehat\zeta_X$ denotes the complete 
Artin zeta function of $X$;

\noindent
(ii) The period of $(G,P)$ for $X$ is defined by
$$\begin{aligned}\omega^{(G,P)}_X(s):=&\mathrm{Res}_{\langle\lambda-\rho,\beta_{P,n-1}^\vee\rangle=0}\cdots
 \mathrm{Res}_{\langle\lambda-\rho,\beta_{P,2}^\vee\rangle=0}\mathrm{Res}_{\langle\lambda-\rho,\beta_{P,1}^\vee\rangle=0}\,\omega_X^G(\lambda)\\
 =&\mathrm{Res}_{s_n=0}\cdots \mathrm{Res}_{s_{p+1}=0}\mathrm{Res}_{s_{p-1}=0}\cdots
 \mathrm{Res}_{s_1=0}\,\omega_X^G(\lambda)\end{aligned}$$ with $s=s_p$.
}
 
 It is clear that there exists a minimal number $I(G,P)$ and factors $$\widehat\zeta_X\Big(a_1^{(G,P)}s+b_1^{(G,P)}\Big),\ \widehat\zeta_X\Big(a_2^{(G,P)}s+b_2^{(G,P)}\Big),\ \cdots,\ 
\widehat\zeta_X\Big(a_{I(G,P)}^{(G,P)}s+b_{I(G,P)}^{(G,P)}\Big),$$
such that there are no zeta factors appeared in the denominators of
all terms of the product $\Big[\prod_{i=1}^{I(G,P)}\widehat\zeta_X\Big(a_i^{(G,P)}s+b_i^{(G,P)}\Big)\Big]\cdot \omega_{\mathbb Q}^{(G,P)}(s).$ 
 
\vskip 0.30cm
\noindent
{\bf Main Definition 3.}  {\it  The zeta function of $X$ associated to $(G,P)$ 
is defined by $${\widehat\zeta}_X^{(G,P)}(s):=\Big[\prod_{i=1}^{I(G,P)}\widehat\zeta_X\Big(a_i^{(G,P)}s+b_i^{(G,P)}\Big)\Big]\cdot \omega_{X}^{(G,P)}(s).$$}

\begin{thm} (Functional Equation) We have
$${\widehat\zeta}_X^{(G,P)}(-c_p-s)={\widehat\zeta}_X^{(G,P)}(s).$$
\end{thm}

\subsection{Proof of the Functional Equation}
Using the Lie structures exposed, next,  we give a proof of the functional equation for the group zetas of function fields, 
following [Ko], in which the group zetas for the field $\mathbb Q$ of rational numbers is treated.
\subsubsection{Lie Structures}
For $w\in W$, denote by $l(w):=|\Phi_w|$ the length of $w$. Write   the longest element of $W$ as $w_0$. Then, $$w_0^2=id,\quad w_0\Delta=\Delta\qquad\mathrm{and}\quad w_0\Phi^+=\Phi^-.$$ Similarly, for a fixed $p$, denote by $w_p$ the longest element of $W_p$.
Now, for $w\in W$, introduce the subset $\W_p$ of $W$
by $$\W_p:=\{w\in W:w\Delta_p\subset\Delta\cup\Phi^-\}.$$ One checks that $\id,w_0,w_p\in\W_p$.

For each $\alpha\in \Phi$, define its height by $\het\,\alpha^\vee:=\langle\rho,\alpha^\vee\rangle$.
For $w\in\W_p$ and $(k,h)\in\mathbb Z^2$, set
$$\begin{aligned}N_{p,w}(k,h):=&\#\{\alpha\in w^{-1}\Phi^-:\langle\lambda_p,\alpha^\vee\rangle=k,\het\,\alpha^\vee=h\},\\
N_{p}(k,h):=&\#\{\alpha\in \Phi:\langle\lambda_p,\alpha^\vee\rangle=k,\het\,\alpha^\vee=h\},\\
M_p(k,h):=&\max_{w\in \W_p}\{N_{p,w}(k,h-1)-N_{p,w}(k,h)\},\qquad\mathrm{and}\\
\widetilde M_p(k,h):=&\max_{w\in \W_p}\{\delta(N_{p,w}(k,h-1)-N_{p,w}(k,h))\},\end{aligned}$$
where $\delta(a)=a$ if $a>0$ and 0 otherwise.

\begin{lem}  The following relations hold.

\noindent
(i)  If $h\geq 1$, $M_p(k,h)=\widetilde M_p(k,h)$;

\noindent
(ii) $N_p(k,kc_p-h)-M_p(k,kc_p-h+1)=N_p(k,h-1)-M_p(k,h)$;

\noindent
(iii)  $c_p\lambda_p-w_p\rho=\rho$.
\end{lem}

They are various lemmas of [Ko]. More precisely, (i), (ii) and (iii) correspond to
Lem. 5.4 (1), (2) and Lem 4.1, respectively.
\subsubsection{A local decomposition}
Write by $$\omega_{X}^G(\lambda)=:\sum_{w\in W}\omega_w^G(\lambda)$$ where
 $$\omega_w^G(\lambda):=\Big(\prod_{\alpha\in\Delta}\frac{1}{1-q^{-\langle w\lambda-\rho,\alpha^\vee\rangle}}\Big)\Big(\prod_{\alpha\in\Phi_w}\frac{\widehat\zeta_X(\langle\lambda,\alpha^\vee\rangle)}
{\widehat\zeta_X(\langle\lambda,\alpha^\vee\rangle+1)}\Big).$$
Since $$\langle w\lambda,\lambda'\rangle=\langle \lambda,w^{-1}\lambda'\rangle,\quad \ w\alpha^\vee=(w\alpha)^\vee,$$ we have, for each $w\in W$, locally, 
$$\begin{aligned}\omega_w^G&(\lambda)=
\Big(\prod_{\alpha\in\Delta}\frac{1}{1-q^{-\langle w\lambda-\rho,\alpha^\vee\rangle}}\Big)
\Big[\prod_{\alpha\in\Phi_w\,\cap\, \Delta_p}\frac{1}{1-q^{-\langle \lambda-\rho,\alpha^\vee\rangle}}
\Big]\\
\times&\Big(\Big[\prod_{\alpha\in\Phi_w\,\cap\, \Delta_p}\big(1-q^{-\langle \lambda-\rho,\alpha^\vee\rangle}\big)\cdot\frac{\widehat\zeta_X(\langle\lambda,\alpha^\vee\rangle)}
{\widehat\zeta_X(\langle\lambda,\alpha^\vee\rangle+1)}\Big]\Big[\prod_{\alpha\in\Phi_w\backslash\Delta_p}\frac{\widehat\zeta_X(\langle\lambda,\alpha^\vee\rangle)}
{\widehat\zeta_X(\langle\lambda,\alpha^\vee\rangle+1)}\Big]\Big)\\
=&\Big[\Big(\prod_{\alpha\in (w^{-1}\Delta\cup \Phi_w)\,\cap\,\Delta_p}\frac{1}{1-q^{-\langle \lambda-\rho,\alpha^\vee\rangle}}\Big)
\Big(\prod_{\alpha\in (w^{-1}\Delta)\backslash \Delta_p}\frac{1}{1-q^{-\langle \lambda-\rho,\alpha^\vee\rangle}}
\Big)\Big]\\
\times&\Big(\Big[\prod_{\alpha\in\Phi_w\,\cap\, \Delta_p}\big(1-q^{-\langle \lambda-\rho,\alpha^\vee\rangle}\big)\cdot\frac{\widehat\zeta_X(\langle\lambda,\alpha^\vee\rangle)}
{\widehat\zeta_X(\langle\lambda,\alpha^\vee\rangle+1)}\Big]\Big[\prod_{\alpha\in\Phi_w\backslash\Delta_p}\frac{\widehat\zeta_X(\langle\lambda,\alpha^\vee\rangle)}
{\widehat\zeta_X(\langle\lambda,\alpha^\vee\rangle+1)}\Big]\Big).\end{aligned}
$$  
\subsubsection{Taking residues}
Next, for $\omega_w^G(\lambda)$, we take the residues at $s_k=0$ for $k\not=p$ and put $s_p=s$. Recall that  $$\Big[\alpha\in\Delta\quad\Leftrightarrow\quad\langle\rho,\alpha^\vee\rangle=1\Big]\qquad\Rightarrow\qquad\langle\lambda-\rho,\alpha^\vee\rangle=\sum_{k=1}^na_ks_k.$$
Consequently, for each of four products appeared in the latest expression for $\omega_w^G(\lambda)$, (after taking the residue), we have

(i) For the first term, $$\prod_{\alpha\in (w^{-1}\Delta\cup \Phi_w)\,\cap\,\Delta_p}\frac{1}{1-q^{-\langle \lambda-\rho,\alpha^\vee\rangle}}=\prod_{\alpha\in (w^{-1}\Delta\cup \Phi_w)\,\cap\,\Delta_p}\frac{1}{1-q^{-s_k}};$$

(ii) For the second term,
$$\prod_{\alpha\in (w^{-1}\Delta)\backslash \Delta_p}\frac{1}{1-q^{-\langle \lambda-\rho,\alpha^\vee\rangle}}
\Big|_{s_k=0, k\not=p;s_p=s}=\prod_{\alpha\in (w^{-1}\Delta)\backslash \Delta_p}\frac{1}{1-q^{-\langle \lambda_p,\alpha^\vee\rangle s-\het\,\alpha^\vee+1}}.$$ Since 
$\alpha\in (w^{-1}\Delta)\backslash \Delta_p$, $\het\,\alpha^\vee\not=1$ or $\langle \lambda_p,\alpha^\vee\rangle\not=0$. Thus the denominator do not vanish identically.

(iii) In the third term, for $\alpha_k\in\Phi_w\,\cap\, \Delta_p$, we have
$$\Big(1-q^{-\langle \lambda-\rho,\alpha_k^\vee\rangle}\Big)\cdot\frac{\widehat\zeta_X(\langle\lambda,\alpha_k^\vee\rangle)}
{\widehat\zeta_X(\langle\lambda,\alpha_k^\vee\rangle+1)}=\big(1-q^{-s_k}\Big)\cdot\frac{\widehat\zeta_X(s_k+1)}{\widehat\zeta_X(s_k+2)}=\frac{\widehat\zeta_X(1)}{\widehat\zeta_X(2)}+o(s_k)$$ as $s_k\to 0$, where
$\widehat\zeta_X(1):=\mathrm{Res}_{s=1}\widehat\zeta_X(s)$.

(iv) In the forth term, for $
\alpha\in\Phi_w\backslash\Delta_p$, we have
$$\frac{\widehat\zeta_X(\langle\lambda,\alpha^\vee\rangle)}
{\widehat\zeta_X(\langle\lambda,\alpha^\vee\rangle+1)}\Big|_{s_k=0,k\not=p;s_p=s}=
\frac{\widehat\zeta_X(\langle\lambda_p,\alpha^\vee\rangle s+\het\,\alpha^\vee)}
{\widehat\zeta_X(\langle\lambda_p,\alpha^\vee\rangle s+\het\,\alpha^\vee+1)}.$$

Consequently, when taking the residues, all terms $\omega_w^G(\lambda)$ vanish
except for the $w$'s satisfying $\Delta_p\subset w^{-1}\Delta\cup \Phi_w$, i.e., $w\in\frak W_p$. Moreover, for $w\in\frak W_p$,
$$\begin{aligned}\Res_{s_k=0,k\not=p}\omega_w^G(\lambda)=&\Big(\prod_{\alpha\in (w^{-1}\Delta)\backslash \Delta_p}\frac{1}{1-q^{-\langle \lambda_p,\alpha^\vee\rangle s-\het\,\alpha^\vee+1}}\Big)\\
&\times\Big(\prod_{\alpha\in\Phi_w\,\cap\, \Delta_p}
\frac{\widehat\zeta_X(1)}{\widehat\zeta_X(2)}\Big)\Big(\prod_{\alpha\in\Phi_w\backslash\Delta_p}\frac{\widehat\zeta_X(\langle\lambda_p,\alpha^\vee\rangle s+\het\,\alpha^\vee)}
{\widehat\zeta_X(\langle\lambda_p,\alpha^\vee\rangle s+\het\,\alpha^\vee+1)}\Big).\end{aligned}$$
Therefore,
$$\begin{aligned}\omega_X^{(G,P)}(s)=&\sum_{w\in \frak W_p}\Big(
\prod_{\alpha\in (w^{-1}\Delta)\backslash \Delta_p}\frac{1}{1-q^{-\langle \lambda_p,\alpha^\vee\rangle s-\het\,\alpha^\vee+1}}\Big)\\
&\times\Big(\prod_{\alpha\in\Phi_w\,\cap\, \Delta_p}
\frac{\widehat\zeta_X(1)}{\widehat\zeta_X(2)}\Big)\Big(\prod_{\alpha\in\Phi_w\backslash\Delta_p}\frac{\widehat\zeta_X(\langle\lambda_p,\alpha^\vee\rangle s+\het\,\alpha^\vee)}
{\widehat\zeta_X(\langle\lambda_p,\alpha^\vee\rangle s+\het\,\alpha^\vee+1)}\Big)\\
&\hskip 1.0cm(\mathrm{since}\quad\Delta_p\subset w^{-1}\Delta\cup \Phi_w\ \Leftrightarrow\ \Delta_p\subset w^{-1}(\Delta\cup \Phi^-))\\
=&\sum_{w\in \frak W_p}\Big(
\prod_{\alpha\in (w^{-1}\Delta)\backslash \Delta_p}\frac{1}{1-q^{-\langle \lambda_p,\alpha^\vee\rangle s-\het\,\alpha^\vee+1}}\Big)\\
&\times\Big(\prod_{\alpha\in\Phi_w\backslash \Delta_p}
\widehat\zeta_X(\langle\lambda_p,\alpha^\vee\rangle s+\het\,\alpha^\vee)\Big)
\Big(\prod_{\alpha\in\Phi_w}\frac{1}
{\widehat\zeta_X(\langle\lambda_p,\alpha^\vee\rangle s+\het\,\alpha^\vee+1)}\Big)\\
=&\sum_{w\in \frak W_p}\Big(
\prod_{\alpha\in (w^{-1}\Delta)\backslash \Delta_p}\frac{1}{1-q^{-\langle \lambda_p,\alpha^\vee\rangle s-\het\,\alpha^\vee+1}}\Big)\cdot \widehat\zeta_{p,w}(s)\\
\end{aligned}\eqno(2)
$$
 where, for $w\in\frak W_p$, we let
 $$\begin{aligned}\widehat\zeta_{p,w}(s):=&\Big(\prod_{\alpha\in\Phi_w\backslash \Delta_p}
\widehat\zeta_X(\langle\lambda_p,\alpha^\vee\rangle s+\het\,\alpha^\vee)\Big)\\
&\times 
\Big(\prod_{\alpha\in\Phi_w}\frac{1}
{\widehat\zeta_X(\langle\lambda_p,\alpha^\vee\rangle s+\het\,\alpha^\vee+1)}\Big).\end{aligned}$$

\subsubsection{Minimal number of factors}
With the above decomposition, we are ready to find out the minimal number of factors used in the normalization process appeared in Main Definition 3. With the expression for $\omega_X^{(G,P)}(s)$ in (2), we concentrate the zeta factors in $\widehat\zeta_{p,w}(s)$ for $w\in\frak W_p$.

By definition,
$$\begin{aligned}\prod_{\alpha\in\Phi_w\backslash \Delta_p}&
\widehat\zeta_X(\langle\lambda_p,\alpha^\vee\rangle s+\het\,\alpha^\vee)\\
=&\widehat\zeta_X(s+1)^{N_{p,w}(1,1)}\prod_{\alpha\in\Phi_w\backslash \Delta}
\widehat\zeta_X(\langle\lambda_p,\alpha^\vee\rangle s+\het\,\alpha^\vee)\\
=&\widehat\zeta_X(s+1)^{N_{p,w}(1,1)}\prod_{k=0}^\infty\prod_{h=2}^\infty \widehat\zeta_X(k s+h)^{N_{p,w}(k,h)},\end{aligned}$$ and
$$\begin{aligned}
\prod_{\alpha\in\Phi_w}&\frac{1}
{\widehat\zeta_X(\langle\lambda_p,\alpha^\vee\rangle s+\het\,\alpha^\vee+1)}\\
=&\prod_{k=0}^\infty\prod_{h=1}^\infty \widehat\zeta_X(k s+h+1)^{-N_{p,w}(k,h)}\\
=&\prod_{k=0}^\infty\prod_{h=1}^\infty \widehat\zeta_X(k s+h+1)^{-N_{p,w}(k,h-1)}.\end{aligned}$$
Hence
$$\widehat\zeta_{p,w}(s)=\widehat\zeta_X(s+1)^{N_{p,w}(1,1)}\prod_{k=0}^\infty\prod_{h=2}^\infty \widehat\zeta_X(k s+h)^{N_{p,w}(k,h)-N_{p,w}(k,h-1)}.$$ 
Therefore,

\hskip 2.0cm$\widehat\zeta_X(k s+h)$ {\it appears in the denominator of} $\omega_X^{(G,P)}(s)$\\
\hskip 6.0cm$\Updownarrow$\\
\hskip 4.0cm $N_{p,w}(k,h)-N_{p,w}(k,h-1)<0$.

Consequently, $$\prod_{k=0}^\infty\prod_{h=2}^\infty \widehat\zeta_X(k s+h)^{\widetilde M_{p}(k,h)}=\prod_{i=1}^{I(G,P)}\widehat\zeta_X\Big(a_i^{(G,P)}s+b_i^{(G,P)}\Big)$$ 
is exactly the minimal zeta factors appeared in the normalization process in defining 
$\widehat\zeta_X^{(G,P)}(s)$.
Thus, by Lem.\,5(i), we have proved  the following

\begin{thm}    The zeta function for $X$ associated to $(G,P)$ is given by
$$\widehat \zeta_X^{(G,P)}(s)=\omega_X^{(G,P)}(s)\cdot \prod_{k=0}^\infty\prod_{h=2}^\infty \widehat\zeta_X(k s+h)^{M_{p}(k,h)}.\eqno(3)$$
\end{thm}
\subsubsection{A global decomposition}
The factor
$$\prod_{k=0}^\infty\prod_{h=2}^\infty \widehat\zeta_X(k s+h)^{M_{p}(k,h)}$$
appeared in (3) proves to be  a bit hard. To overcome this,  we go back to
the expression of $\omega_X^{(G,P)}(s)$ in (2). Introduce the \lq overdone' maximal factor
 $$\begin{aligned}M_X^{(G,P)}(s):=M_p(s):=&\prod_{\alpha\in\Phi^+}\widehat\zeta_X(\langle\lambda_p,\alpha^\vee\rangle s+\het\,\alpha^\vee+1)\\
 =&\prod_{\alpha\in\Phi^-}\widehat\zeta_X(\langle\lambda_p,\alpha^\vee\rangle s+\het\,\alpha^\vee)\end{aligned}.$$ Obviously, 
 being maximal, $M_p(s)$ does clear up all the zeta factors in the denominators of terms of (2).
Moreover, by definition, 
 $$M_p(s)=\prod_{k=0}^\infty\prod_{h=1}^\infty \widehat\zeta_X(k s+h+1)^{N_{p}(k,h)}
 =\prod_{k=0}^\infty\prod_{h=2}^\infty \widehat\zeta_X(k s+h)^{N_{p}(k,h-1)}.$$ Here, in the last step, we have used the functional equation for Artin zetas. This then leads to the global decomposition
 $$\widehat \zeta_X^{(G,P)}(s)=\frac{\Omega_X^{(G,P)}(s)}{D_X^{(G,P)}(s)}$$
where we have set
$$\begin{aligned}\Omega_X^{(G,P)}(s):=&M_X^{(G,P)}(s)\cdot \omega_X^{(G,P)}(s),\quad\mathrm{and}\\
D_X^{(G,P)}(s):=& \prod_{k=0}^\infty\prod_{h=2}^\infty \widehat\zeta_X(k s+h)^{-M_{p}(k,h)+N_{p,w}(k,h-1)}.\end{aligned}$$

As such, then the functional equation of our zeta functions is equivalent to the following

\begin{prop}
$$D_X^{(G,P)}(-c_P-s)=D_X^{(G,P)}(s),\qquad\mathrm{and}\qquad \Omega_X^{(G,P)}(-c_P-s)=\Omega_X^{(G,P)}(s).$$
\end{prop}

\subsubsection{Functional Equation for $D_X^{(G,P)}(s)$}
This is rather easy.
Decompose $D$ according to whether it consists of special values of zetas or not to get $$D_X^{(G,P)}(s):=D_p^0\cdot
D_p^1(s)$$ where
$$\begin{aligned}D^0_p:=&\prod_{h=2}^\infty\widehat\zeta_X(h)^{N_p(0,h-1)-M_p(0,h)},\\
D_p^1(s):=&\prod_{k=1}^\infty\prod_{h=2}^\infty\widehat\zeta_X(ks+h)^{N_p(k,h-1)-M_p(k,h)}.\end{aligned}$$
It suffices to show that 
$$D_p^1(-c_p-s)=D_p^1(s).$$
Since $N_{p,w}(k,h-1)=0$ and $M_p(k,h)=0$ for $k\geq 1$ and $h\leq 1$, we have
$$D_p^1(s)=\prod_{k=1}^\infty\prod_{h=-\infty}^\infty\widehat\zeta_X(ks+h)^{N_p(k,h-1)-M_p(k,h)}.$$ Consequently,
$$\begin{aligned}D_p^1(-c_p-s)=&\prod_{k=1}^\infty\prod_{h=-\infty}^\infty\widehat\zeta_X(-kc_p-ks+h)^{N_p(k,h-1)-M_p(k,h)}\\
=&\prod_{k=1}^\infty\prod_{h=-\infty}^\infty\widehat\zeta_X(ks+kc_p-h+1)^{N_p(k,h-1)-M_p(k,h)}\\
&\qquad(\mathrm{by\ the\ functional\ equation}\ \widehat\zeta_X(1-s)=\widehat\zeta_X(s))\\
=&\prod_{k=1}^\infty\prod_{h=-\infty}^\infty\widehat\zeta_X(ks+h)^{N_p(k,kc_p-h)-M_p(k,kc_p-h+1)}\\
=&\prod_{k=1}^\infty\prod_{h=-\infty}^\infty\widehat\zeta_X(ks+h)^{N_p(k,h-1)-M_p(k,h)}\ (\mathrm{by\ Lem\,5(ii)})\\
=&D_p^1(s).\end{aligned}$$

\subsubsection{Involution Structure on $\frak W_p$}

We are left with the proof of the functional equation for $\Omega_X^{(G,P)}(s)$. For this, 
we use an involution structure on $\frak W_p$ given by $w\mapsto w_0ww_p$.
Set  then $$\begin{aligned}
f_{p,w}(s):=&\prod_{\alpha\in(w^{-1}\Delta\backslash\Delta_p)}\frac{1}{1-q^{-\langle\lambda_p,\alpha^\vee\rangle s-\het\,\alpha^\vee+1}},\\
g_{p,w}(s):=&\prod_{\alpha\in(w^{-1}\Phi^-)\backslash\Delta_p}\widehat\zeta_X(\langle\lambda_p,\alpha^\vee\rangle s+\het\,\alpha^\vee).\end{aligned}$$
Similarly as in [Ko], we have the following
\begin{prop}
(i) ({\bf Involution Structure})
$$f_{p,w}(-c_p-s)=f_{p,w_0ww_p}(s),\qquad g_{p,w}(-c_p-s)=g_{p,w_0ww_p}(s);$$ 

\noindent
 (ii)   $$
\Omega_X^{(G,P)}(s)=
\sum_{w\in \frak W_p} f_{p,w}(s)\cdot g_{p,w}(s).$$
\end{prop}

\noindent
{\it Proof.} (i)  For a fixed subset $A\subset\Phi, w\in W$, set
$$S_{p,A}(s;w):=\big\{\langle\lambda_p,\alpha^\vee\rangle s+\het\,\alpha^\vee:\alpha\in(w^{-1}A)\backslash \Delta_p\big\}.$$
Then, note that, for $A=\Delta$ or $\Phi^-$, 
$w_0A=-A$ and 
$$-w_p(w^{-1}A\backslash\Delta_p)=(w_pw^{-1}(-A))\backslash (w_p(-\Delta_p))=(w_pw^{-1}w_0A)\backslash \Delta_p.$$ So, we have $$\begin{aligned}f_{p,w}(s)=&\prod_{as+b\in S_{p,\Delta}(s;w)}\frac{1}{1-q^{-as-b+1}},\\
g_{p,w}(s)=&\prod_{as+b\in S_{p,\Phi^-}(s;w)}\frac{1}{\widehat\zeta_X(as+b)},\end{aligned}$$
 Moreover, 
$$\begin{aligned}S_{p,A}(-c_p-s;w)=&\{\langle\lambda_p,\alpha^\vee\rangle(-c_p-s)+\het\,\alpha^\vee:\alpha\in (w^{-1}A)\backslash \Delta_p\}\\
=&\{\langle\lambda_p,-w_p\alpha^\vee\rangle s+\langle c_p\lambda_p-w_p\rho,-w_p\alpha^\vee\rangle
:\alpha\in (w^{-1}A)\backslash \Delta_p\}\\
=&\{\langle\lambda_p,\beta^\vee\rangle s+\langle \rho,\beta^\vee\rangle
:\beta\in (w_pw^{-1}w_0A)\backslash \Delta_p\}\\
&\hskip 4.0cm (\mathrm{by\ Lem\, 5(iii)})\\
=&S_{p,A}(s;w_0ww_p).\end{aligned}$$ 

\noindent
(ii) In [Ko], the following Lie structures are exposed.
 $$\begin{aligned}(a)\ \Phi^-\backslash (-\Phi_w)&=\Phi^-\backslash (\Phi^-\,\cap\, w^{-1}\Phi^+)=\Phi^-\backslash w^{-1}\Phi^+=\Phi^-\,\cap\, w^{-1} \Phi^-,\\
(b)\ \ (\Phi_w\backslash\Delta_p)\cup&(\Phi^-\,\cap\, w^{-1}\Phi^-)=
 ((\Phi^+\cap w^{-1}\Phi^-)\backslash\Delta_p)\cup(\Phi^-\,\cap\, w^{-1}\Phi^-)\\
 =& ((\Phi^+\cap w^{-1}\Phi^-)\cup(\Phi^-\,\cap\, w^{-1}\Phi^-))\backslash\Delta_p
 =w^{-1}\Phi^-\,\backslash\,\Delta_p.\end{aligned}$$
Consequently,
$$\begin{aligned}\Omega_X^{(G,P)}(s)\buildrel {(a)}\over=&\sum_{w\in W,\Delta_p\subset w^{-1}(\Delta\cup\Phi^-)}\Big(\prod_{\alpha\in(w^{-1}\Delta)\backslash\Delta_p}\frac{1}{1-q^{-\langle\lambda_p,\alpha^\vee\rangle s-\het\,\alpha^\vee+1}}\Big)\\
&\hskip 1.5cm\times\Big(\prod_{\alpha\in(\Phi_w\backslash\Delta_p)\cup(\Phi^-\,\cap\, w^{-1}\Phi^-)}\widehat\zeta_X\big(\langle\lambda_p,\alpha^\vee\rangle s+\het\,\alpha^\vee\big)\Big)\\
\buildrel {(b)}\over=&\sum_{w\in \frak W_p} \Big(\prod_{\alpha\in(w^{-1}\Delta\backslash\Delta_p)}\frac{1}{1-q^{-\langle\lambda_p,\alpha^\vee\rangle s-\het\,\alpha^\vee+1}}\Big)\\
&\hskip 1.5cm\times\Big(\prod_{\alpha\in(w^{-1}\Phi^-)\backslash\Delta_p}\widehat\zeta_X(\langle\lambda_p,\alpha^\vee\rangle s+\het\,\alpha^\vee)\Big).\end{aligned}$$

\noindent
{\it Proof for Thm 4.}
With Prop. 8, the functional equation $$\Omega_X^{(G,P)}(c_P-s)=\Omega_X^{(G,P)}(s)$$ 
is a direct consequence of the fact that $w_0,w_p\in\frak W_p$ so $w\mapsto w_0ww_p$ induces an involution structure of $\frak W_p$. This then completes the proof of Thm 4 as well.
\vskip 0.30cm
\noindent
{\it Remark.} The  functional equation for $\widehat\zeta_X^{(G,P)}(s)$ can be more directly proved using the involution structure $w\mapsto w_0w w_p$ on $\W_p$
directly, based on  the relation ([KKS, Cor 8.7])
$$M_p(k,h)=N_{p,w_0}(k,h-1)-N_{p,w_0}(k,h).$$

\section{Counting Bundles}
\subsection{Uniformity and the Riemann Hypothesis}
Recall that, by Thm 2, or better, the equality (1), we have
$$\begin{aligned}\zeta_{X,r}&(s)
=\sum_{m=0}^{(g-1)-1}\alpha_{X,r}(mr)
\cdot \Big((q^{-rs})^{m}+ (q^r)^{(g-1)-m}\cdot (q^{-rs})^{2(g-1)-m}\Big)\\
&+\alpha_{X,r}\big(r(g-1)\big)\cdot (q^{-rs})^{g-1}+(q^r-1)\beta_{X,r}(0)\cdot \frac{(q^{-rs})^{g}}{(1-q^{-rs})(1-q^rq^{-rs})}.
\end{aligned}$$
Thus, the following conjecture, motivated by our works on zetas for number fields ([W2,3]), 
counts semi-stable bundles decisively.

\begin{conj} (1) ({\bf Uniformity}) There are universal constants $a_{F,r}, b_{F,r}$ and rational functions  $c_{F,r}(q)$ depending on $F$ and $r$ such that
$$\widehat\zeta_{F,r}(s)=c_{F,r}(q)\cdot \widehat\zeta_X^{(SL_r,P_{r-1,1})}(a_{F,r}\cdot s+b_{F,r}).$$

\noindent
(2) ({\bf The Riemann Hypothesis})
$${\widehat\zeta}_{F,r}(s)=0\qquad\Rightarrow\qquad\mathrm{Re}(s)=\frac{1}{2}.$$
\end{conj}

That is to say,  all weighted counts on semi-stables via the invariants $\alpha$'s and $\beta$'s can be read from Artin's zeta functions defined using only line bundles, while  the Riemann Hypothesis 
gives an effective control of the invariants $\alpha$'s and $\beta$'s.

We have the following supportive evidences.

\begin{thm} (i) ({\bf Uniformity}, [W4]) For elliptic curves, the uniformity holds when $r=1,\,2,\,3,\,4,\,5.$

\noindent
(ii) ({\bf Riemann Hypothesis})  The Riemann Hypothesis holds for

\noindent
(a) (Weil) $\widehat\zeta_X(s)$;

\noindent
(b) ([Y]) $\widehat\zeta_X^{(SL_2,P_{1,1})}(s)$;

\noindent
(c) ([W4]) $\widehat\zeta_{E,r}(s)$ for $r=2,3,4,5$ with $E$ an elliptic curve.
\end{thm}

\subsection{Parabolic Reduction, Stability and the Mass}
To end this paper, we explain the reasons why $\zeta_{X,r}(s)$ are {\it non-abelian} zeta functions of $X$, despite the uniformity claiming that, up to certain rational function factors,
$\zeta_{X,r}(s)$ can be read from {\it abelian} Artin zetas.

The central reason is certainly that $\zeta_{X,r}(s)$'s  are defined using moduli spaces of semi-stable bundles, highly non-commutative objects associated to $X$. Furthermore, even assuming the uniformity, from the equation (2), we can still detect where the non-abelian structure lies on. More precisely, in each term $\omega_w^{G}(s),\ w\in\frak W_p$, while, for the zeta factor part
$$\begin{aligned}\widehat\zeta_{p,w}(s):=&\Big(\prod_{\alpha\in\Phi_w\backslash \Delta_p}
\widehat\zeta_X(\langle\lambda_p,\alpha^\vee\rangle s+\het\,\alpha^\vee)\Big)\\
&\times 
\Big(\prod_{\alpha\in\Phi_w}\frac{1}
{\widehat\zeta_X(\langle\lambda_p,\alpha^\vee\rangle s+\het\,\alpha^\vee+1)}\Big),\end{aligned}$$ there are involved only Artin zetas which are abelian,  the non-abelian structure, through the group structure,
is naturally reflected via the
rational function factors $$\prod_{\alpha\in (w^{-1}\Delta)\backslash \Delta_p}\frac{1}{1-q^{-\langle \lambda_p,\alpha^\vee\rangle s-\het\,\alpha^\vee+1}}.$$

To properly understand this, let us examine the so-called parabolic reduction structure appeared in the mass formula for function fields, and similarly, the volume formula for number fields, respectively.

As usual, let $\mathcal D_{\mathbb Q,r}$ be the volume of fundamental domain of $SL_r(\mathbb Z)$, and $\mathcal M_{\mathbb Q,r}[1]$ the moduli space of rank $r$
semi-stable lattices of volume 1. (For background materials, please refer to [W2].)
Then we have

\begin{thm} (i) (Siegel) $$\mathrm{Vol}\Big(\mathcal D_{\mathbb Q,r}\Big)=r\cdot\prod_{i=1}^r\widehat\zeta(i);$$

\noindent
(ii) 
(Reformulation of [KS, \S16]) 
$$\mathrm{Vol}\Big(\mathcal D_{\mathbb Q,r}\Big)=\sum_{\substack{n_1,\dots,n_k\geq 1,\\ n_1+\cdots+n_k=r}}\frac{\prod_{j=1}^k\mathrm{Vol}\Big(\mathcal M_{\mathbb Q,n_j}[1]\Big)}{n_1(n_1+n_2)\cdots(n_1+\cdots+n_k)\cdots(n_k-1+n_k)n_k};$$ 

\noindent
(iii) (Weng [W2, \S4.8]) $$\frac{1}{r}\cdot \mathrm{Vol}\Big(\mathcal M_{\mathbb Q,r}[1]\Big)=\sum_{\substack{n_1,\dots,n_s>0,\\ n_1+\cdots+n_k=r}}
\frac{(-1)^{k-1}}{\prod_{j=1}^{k-1}(n_j+n_{j+1})}\prod_{j=1}^k \mathrm{Vol}\Big(\mathcal D_{\mathbb Q,n_j}\Big).$$
\end{thm}

\noindent
{\it Remarks.} (1) Siegel's formula claims that the volume of {\it non-abelian} fundamental domain can be
measured using special values of the {\it abelian} zeta;

\noindent
(2) Even the roots for [KS, \S16] and [W2] are very much different: the former uses arithmetic truncation of Harder-Narasimhan filtration, and the later
uses analytic truncation and Eisenstein series, they share a common origin, as we observed, namely,
the {\it parabolic reduction structure};

\noindent
(3) The part of non-abelian group structure and the part of the abelian zeta  are well-organized so that they fit into a uniform theory naturally. For example, roughly, we see that fundamental domains consists of an essential part coming from stable lattices and boundary parts coming from tubular  neighborhoods of cusps associated to proper parabolic subgroups.

Motivated by this, more generally, for  a split reductive group $G$  defined over a number field $F$, $B$ a fixed Borel ...
denote by $G(\mathbb A)^{\mathrm{ss}}$ the adelic elements of $G$ corresponding to semi-stable principle $G$-lattices ([G]). Write $\mathbb K_G$ for the associated maximal compact subgroup.
 Also for a standard parabolic subgroup $P$, write its  Levi decomposition as  $P=UM$ with $U$ the unipotent radical and $M$ its Levi factor. Denote the corresponding simple decomposition of $M$ as
$\prod_iM_i$ with $M_i$'s the simple factors of $M$.
Introduce invariants $$\nu_P:=\prod_i\Vol\Big(\mathbb K_{M_i} Z_{M_i^1(\mathbb A)}\big\backslash M_i^1(\mathbb A)\big/M_i(F)\Big)$$
and $$\mu_P:=\prod_i\Vol\Big(\mathbb K_{M_i} Z_{M_i^1(\mathbb A)}\big\backslash M_i^1(\mathbb A)^{\mathrm{ss}}\big/M_i(F)\Big).$$ In parallel, we have similar constructions for function fields $F=\mathbb F_q(X)$. Based on all this, then we have the following

\begin{conj} ({\bf Parabolic Reduction}) Let $G/F$ be a split reductive group with $B/F$ a fixed Borel. Then, for each standard parabolic subgroup $P$ of $G$, there exist constants 
$c_P\in\mathbb Q,\ e_P\in\mathbb Q_{>0}$ such that
 $$\nu_G=\sum_{P}\,c_P\cdot\nu_P,\qquad
 \mu_G=\sum_{P}\,\mathrm{sgn}(P)\,e_P\cdot\nu_P,$$ where $P$ runs over all standard parabolic subgroups of $G$, and $\mathrm{sgn}(P)$ denotes the sign of $P$.  
 \end{conj}
 
The exact values of $e_P$'s can be written out in terms of the associated root system. 
Indeed, 
if $$W_0:=\Big\{w\in W:\{\alpha\in\Delta:w\alpha\in \Delta\cup \Phi^- \}=\Delta\Big\},$$
then there is a natural one-to-one correspondence between $W_0$ and the set of subsets of $\Delta$, and hence to the set of standard parabolic subgroups of $G$. Thus we will write
$$W_0:=\Big\{w_P:\ P\ \mathrm{standard\ parabolic\ subgroup}\Big\},$$
and, for $w=w_P\in W_0$, write  $J_P\subset\Delta$ the corresponding subset. 

\begin{conj} ({\bf Parabolic Reduction, Stability \& the Mass}, [W5]) Let $G$ be a split type reductive group
with $P$ its maximal parabolic subgroup.

\noindent
(1) Over a number field $F$,

\noindent
(i) The volume of moduli space of semi-stable principal lattices is given by 
$$\nu_G=\mathrm{Res}_{s=-c_P}\widehat\zeta_F^{(G,P)}(s)=\mathrm{Res}_{\lambda=\rho}\omega_F^{G}(\lambda);$$

\noindent
(ii) We have the following formula
$$ \mu_G=\sum_{P}\frac{(-1)^{\mathrm{rank}(P)}}{\prod_{\alpha\in\Delta\backslash w_JJ_P}(1-\langle w_J\rho,\alpha^\vee\rangle)}\cdot\nu_P;$$

\noindent
(2) Over an irreducible reduced regular projective curve $X$, 

\noindent
(i) The mass of moduli space of semi-stable principal bundles is given by 
$$\log q\cdot  \nu_G=\mathrm{Res}_{s=-c_P}\widehat\zeta_X^{(G,P)}(s)=\mathrm{Res}_{\lambda=\rho}\omega_X^{G}(\lambda);$$

\noindent
(ii) We have the following formula
$$\mu_G=\sum_{P}\frac{(-1)^{\mathrm{rank}(P)}}{\prod_{\alpha\in\Delta\backslash w_JJ_P}(1-q^{\langle w_J\rho,\alpha^\vee\rangle-1})}\cdot\nu_P;$$ 
 \end{conj}

\noindent
{\it Remarks.} (1) We expect that $c_P>0$ for and only for number fields.

\noindent
(2) Calculations in [Ad] for lower ranks groups indicates that, for number fields, $\frac{1}{c_P}\in\mathbb Z_{>0}$. It would be very interesting to find a close formula for them.
\vskip 0.30cm
All this indicates that non-commutative group structures are naturally embedded into our pure high rank zeta functions.
\vskip 0.70cm
\centerline {\bf\large{ REFERENCES}}
\vskip 0.40cm 
\noindent
[Ad] K. Adachi, Parabolic Reduction, Stability and Volumes of Fundamental Domains: Examples of Low Ranks, manuscript, 2012
\vskip 0.10cm
\noindent
[AB] M. Atiyah \& R. Bott, The Yang-Mills equations over Riemann surfaces, Philos. Trans. Roy. Soc. London 308, 523-615 (1983)
\vskip 0.10cm
\noindent
[DR] U.V. Desale \& S. Ramanan, Poincare polynomials of the variety of stable bundles, Math. Ann 26 (1975) 233-244
\vskip 0.10cm
\noindent
[D] V.G. Drinfeld, Number of two-dimensional irreducible representations of the fundamental group of a curve over a finite field,
Func. Anal. \& App., {\bf 15}  (1981), 294-295
 \vskip 0.10cm
\noindent 
[HN] G. Harder \& M.S. Narasimhan, On the cohomology groups of moduli spaces of vector bundles on curves, Math. Ann. 212, 215-248 (1975)
\vskip 0.10cm
\noindent 
[Hu] J.E. Humphreys, {\it Linear Algebraic Groups}, GTM 21, Springer-Verlag (1975)
 \vskip 0.10cm
\noindent 
[KKS] H. Ki, Y. Komori \& M. Suzuki, On the zeros of Weng zeta functions for Chevalley groups, arXiv:1011.4583v1
\vskip 0.20cm  
\noindent
[Ko] Y. Komori, Functional equations for Weng's zeta functions for $(G,P)/\mathbb{Q}$, Amer. J. Math., to appear
\vskip 0.10cm
\noindent
[KS] M. Kontsevich \& Y. Soibelman, Lectures on motivic Donaldson-Thomas invariants and Wall-crossing formulas, manuscripts, Dec.  2011
\vskip 0.10cm
\noindent
[L] L. Lafforgue, {\it Chtoucas de Drinfeld et conjecture de 
Ramanujan-Petersson}. Asterisque No. 243 (1997)
\vskip 0.10cm
\noindent
[La] R. Langlands, The volume of the fundamental domain for some arithmetical 
subgroups of Chevalley groups, in {\it Algebraic Groups and Discontinuous 
Subgroups,} Proc. Sympos. Pure Math. 9, AMS (1966) pp.143--148
\vskip 0.20cm 
\noindent
[W1] L. Weng, Non-abelian zeta functions for function fields. {\it Amer. J. Math}. 127 (2005), no. 5, 973--1017.
 \vskip 0.20cm 
\noindent 
[W2] L. Weng, A geometric approach to $L$-functions. {\it The Conference on $L$-Functions,} 219--370, World Sci. Publ., Hackensack, NJ, 2007.
 \vskip 0.20cm 
\noindent 
[W3] L. Weng, Symmetries and the Riemann Hypothesis, {\it Advanced Studies in Pure Mathematics 58}, Japan Math. Soc., 173-224, 2010
 \vskip 0.20cm 
\noindent 
[W4] L. Weng, Zeta Functions for Elliptic Curves I: Counting Bundles, arXiv:1202.0870
 \vskip 0.20cm 
\noindent 
[W5] L. Weng, Parabolic Reduction, Stability and The Mass, in preparation
\vskip 0.10cm
\noindent
[Wi] E. Witten, On quantum gauge theories in two dimensions, Commun. Math. Phys, 141 (1991) 153-209
\vskip 0.10cm
\noindent
[Y] H. Yoshida, manuscripts, Dec., 2011
\vskip 0.10cm
\noindent
[Z] D. Zagier, Elementary aspects of the Verlinde formula and the Harder-Narasimhan-Atiyah-Bott formula, in {\it Proceedings of the Hirzebruch 65 Conference on Algebraic Geometry}, 445-462 (1996)
\vskip 0.30cm
\noindent
Lin WENG\footnote{{\bf Acknowledgement.}
We would like to thank IHES for providing us an excellent working environment during our September visit in 2011. Special thanks also due to Kontsevich and Yoshida for sharing with us their works, which motivate
our current works,
to Lafforgue for enlightening us the interchangeability between taking residues and taking integrations over moduli spaces of semi-stable lattices, which leads to the uniformity conjecture, and to Deninger and Hida for constant encouragements.

This work is partially supported by JSPS.}
\noindent
Institute for Fundamental Research, The $L$ Academy {\it and}

\noindent
Graduate School of Mathematics, Kyushu University, Fukuoka, 819-0395, 
JAPAN

\noindent
E-Mail: weng@math.kyushu-u.ac.jp

\end{document}